\documentclass[12pt,centertags]{amsart}
\usepackage{amsmath    }
\usepackage{amstext    }
\usepackage{amsthm    }
\usepackage{a4}
\usepackage{amssymb}
\usepackage{amscd}
\usepackage[mathscr]{eucal} 
\usepackage{mathrsfs}
\numberwithin{equation}{section}
\newcommand{\comment}[1]{} 

\newcommand{\field}[1]{\mathbb{#1}}
    \newcommand{\R}{\field{R}} \newcommand{\B}{\field{B}}
\newcommand{\C}{\field{C}}      \newcommand{\Proj}{\field{P}}
\newcommand{\cali}[1]{\mathscr{#1}} 
\newcommand{\cC}{\cali{C}}  
\newcommand{\cO}{\cali{O}}  \newcommand{\cH}{\cali{H}}  \newcommand{\cF}{\cali{F}} 
   
\newcommand{\iom}{\Omega}   

\DeclareMathOperator{\Dom}{Dom}       \DeclareMathOperator{\Ran}{Range}

\DeclareMathOperator{\ke}{Ker} \DeclareMathOperator{\Id}{Id} 
\DeclareMathOperator{\supp}{supp} 
\DeclareMathOperator{\sing}{Sing} \DeclareMathOperator{\reg}{Reg} 
\newcommand{\norm}[1]{\lVert#1\rVert}
\newcommand{\abs}[1]{\lvert#1\rvert}
\newcommand{\om}{\omega}      
\newcommand{\imat}{\sqrt{-1}} 
\newcommand{\db}{\overline\partial} \newcommand{\dbs}{\overline\partial^{\ast}}
\newcommand{\pa}{\partial}  
\newcommand{\lap}{\Delta''}

\newtheorem{theorem}{Theorem}[section]
\newtheorem{lemma}[theorem]{Lemma}
\newtheorem{proposition}[theorem]{Proposition}
\newtheorem{corollary}[theorem]{Corollary}

\theoremstyle{definition} 
\newtheorem{defn}[theorem]{Definition}
\theoremstyle{remark}
\newtheorem{remark}[theorem]{Remark}

\theoremstyle{remark}

\theoremstyle{remark}

\theoremstyle{remark} 
 
\begin{document}
\title[Compactification of hyperconcave ends]{On the compactification of hyperconcave ends and the theorems of
Siu\,-Yau and Nadel}
\author[G. Marinescu]{George Marinescu}  
\thanks{Partially supported by SFB 288 and the Alexander von Humboldt Stiftung} 
\address{Humboldt-Universit{\"a}t zu Berlin, Institut f{\"u}r Mathematik, 
Unter den Linden 6,10099 Berlin, Germany and
Institute of Mathematics of the Romanian Academy, Bucharest, Romania}
\email{george@mathematik.hu-berlin.de} 
\urladdr{} 
\author[T. C. Dinh]{Tien-Cuong Dinh} 
\address{Universit\'e Paris-Sud, B\^atiment 425-Math\'ematique,
91405 Orsay, France} 
\curraddr{} 
\email{Tiencuong.Dinh@math.u-psud.fr} 
\urladdr{}
\date{\today}
\begin{abstract}
We show that the `pseudoconcave holes' of some naturally 
arising class of manifolds, called hyperconcave ends, 
can be filled in, including the case of complex dimension two. 
As a consequence we obtain a stronger version of the  
compactification theorem of Siu\,-Yau and extend Nadel's 
theorems to dimension two.
\end{abstract} 

\maketitle

\section{Introduction} \label{s1} 

We will be concerned with the following class of manifolds. 
\begin{defn}\label{conc}
A complex manifold $X$ with $\dim X\geqslant2$ is said to be a {\em strongly pseudoconcave end} 
if there exist $a\in\R\cup\{+\infty\}$ and $c\in\R\cup\{-\infty\}$ and a proper, smooth function $\varphi:X\longrightarrow(c,a)$, which is strictly plurisubharmonic  
on a set of the form $\{\varphi<b\}$, for some $b\leqslant a$.
If $c=-\infty$, $X$ is called a {\em hyperconcave end}. 
For $e<d<a$ we set $X_d=\{\varphi<d\}$ and $X^d_e=\{e<\varphi<d\}$.  
We call $\varphi$ {\em exhaustion function}.
\end{defn}
We say that a strongly pseudoconcave end can be {\em compactified} or
{\em filled in} if there exists a  
complex space $\widehat{X}$ such that $X$ is (biholomorphic to) an open set  
in $\widehat{X}$ and for any $d<a$, $(\widehat{X}\smallsetminus X)
\cup\{\varphi\leqslant d\}$ 
is a compact set. 
We will call $\widehat{X}$ the {\em completion} of $X$.
 
By a theorem of Rossi \cite[Th.\,3, p.\,245]{Ros:65} and Andreotti-Siu 
\cite[Prop.\,3.2]{AS:70} any strongly pseudoconcave end $X$ can be 
compactified, provided $\dim{X}\geqslant3$. 
This is no longer true if $\dim{X}=2$, as shown in a counterexample of 
Grauert, Andreotti-Siu and Rossi \cite{Gra:94,AS:70,Ros:65}.

Our goal is to compactify the hyperconcave ends also in dimension two.
Let us mention some examples.
The regular part of a variety with isolated singularities is a  
hyperconcave end. The same is true for the complement of a compact  
completely pluripolar set (the set where a strongly  
plurisubharmonic function equals $-\infty$) in a complex manifold.   
The first step in the proof of the Siu-Yau  
compactification theorem \cite{SY:82} is to show that 
a complete K\"ahler manifold $X$ of finite volume and sectional curvature  
pinched between two negative constants has hyperconcave ends. 
One can also check that the examples of Grauert, Andreotti-Siu and Rossi 
are not hyperconcave ends. 

\begin{theorem}\label{main} 
Any hyperconcave end $X$ can be compactified.
Moreover, if $\varphi$ is strictly plurisubharmonic on 
the whole $X$, the completion $\widehat{X}$ can be chosen a normal Stein space with at 
worst isolated singularities.
\end{theorem}  
The motivation for the study of the compactification of hyperconcave ends  
comes from the theory of complex-analytic compactification of quotients 
$X=\B^n/\Gamma$ of the unit ball in $\C^n$, $n\geqslant2$  
by arithmetic groups $\Gamma$. The  
Satake-Baily-Borel compactification $\widehat{X}$ of $X=\B^n/\Gamma$ 
is obtained by adding a finite set of points which are isolated  
singularities. 
The Siu-Yau theorem gives a differential geometric 
proof of this fact, first by proving that $X$ has hyperconcave ends 
and then showing it can be compactified by adding finitely many points. 
In this context, our next goal is to find necessary conditions for a manifold with  
hyperconcave ends to be analytically compactified by adding one point 
at each end. This also yields a complex analytic proof of the second step 
of Siu-Yau's theorem (see Corollary \ref{SiuYau}).
For more results on the compactification of complete K\"ahler-Einstein manifolds 
of finite volume and bounded curvature we refer to Mok
\cite{Mok:89} and the references therein. 
\begin{theorem} \label{nadel}
Let $X$ be a hyperconcave end and let  
$\widehat{X}$ be a smooth completion of $X$.
Assume that $X$ can be covered by Zariski-open sets which are uniformized by 
Stein manifolds. Then $\widehat{X}\smallsetminus{X}$  
is the union of a finite set $D'$ and an exceptional analytic set which  
can be blown down to a finite set $D$.  
Each connected component of $X_c$, for sufficiently small $c$, can be 
analytically compactified by one point from $D'\cup D$. 
If $X$ itself has a Stein cover, 
$D'=\varnothing$ and $D$ consists of the singular set of the Remmert reduction
of $\widehat{X}$.
\end{theorem} 
Let us note that the result is natural, in the light of a recent result of 
Col\c{t}oiu-Tib\u{a}r \cite{CoTi:02}, asserting that the universal cover of a small 
punctured neighbourhood of an isolated singularity of dimension two is Stein, 
whenever the fundamental group of the link is infinite.

Theorem \ref{nadel} affords the extension to dimension two of 
Nadel's theorems \cite{Na:90}. If in Definition \ref{conc} we may  
take $a=+\infty$ and $\varphi$ bounded from above, 
the manifold $X$ is said to be {\em hyper $1$-\,concave}\,.
In this case $\{\varphi\geqslant c\}$ is compact for every $c\in\R$.
Hyper $1$-\,concave manifolds are particular cases of strongly 1-concave 
manifolds in the sense of Andreotti-Grauert \cite{AG:62}.   
\begin{corollary}\label{nad2}
Let $X$ be a connected manifold of dimension $n\geqslant 2$.  
The following conditions 
are necessary and sufficient for $X$ to be a quasiprojective manifold
which can be compactified to a Moishezon space by adding finitely many  
points:
{\rm(}i\,{\rm)} $X$ is hyper $1$-concave,
{\rm(}ii\,{\rm)} $X$ admits a line bundle $E$ such that the ring $\oplus_{k>0}H^0(X,E^k)$ 
separates points and gives local coordinates, and
{\rm(}iii\,{\rm)} $X$ can be covered by Zariski-open sets which can be uniformized by Stein manifolds.
If $X$ has a Stein cover, one adds only singular points. 
\end{corollary}
Corollary \ref{nad2} yields, in dimension two, a version  
of Nadel-Tsuji theorem \cite{NT:88} 
together with a completely complex-analytic proof
of the compactification of arithmetic quotients.  
It answers also \cite[Prob.\,1]{Mok:89} for the case $q=0$. 

In \cite{MY04} we apply the results of this paper to ball quotients of dimension two
having a strongly pseudoconvex boundary, extending results of Burns and Napier-Ramachandran
\cite{NR}. Another application of the above theorems is the embedding of Sasakian manifolds
\cite{MY04,OV:04}. 

The organization of the paper is as follows. In \S \ref{existence} we construct 
holomorphic functions on a hyperconcave end and then we prove Theorem \ref{main}
in \S 3. The proof of Theorem \ref{nadel} occupies \S 4 and in \S5 
we extend Nadel's and Andreotti-Siu theorems to dimension two.

\section{Existence of holomorphic functions} \label{existence}

The idea of proof is to analytically embed small strips  
$X_{e}^{e^*}$, for $e<e^*$ in a neighbourhood of minus infinity,  
into the difference of two concentric polydiscs in the euclidian space.  Then we apply  
the Hartogs extension theorem to extend the image to an analytic set 
which will provide the compactification. 
To obtain the embedding we follow the strategy of  
Grauert and Kohn for the solution of the Levi problem. 
Namely, we solve the $L^2$ $\db$-Neumann for $(0,1)$-forms on  
domains $X_d$ with strongly pseudoconvex boundary $\{\varphi=d\}$ endowed with a 
complete metric at minus infinity. 
This entails the finiteness of the $L^2$ Dolbeault cohomology  
$H^{0,1}_{(2)}(X_d)$ 
which in turn implies the existence of peak holomorphic functions 
at each point of the boundary $\{\varphi=d\}$. 
Note that, as a consequence of the Andreotti-Grauert theory \cite{AG:62}, 
in dimension two, the sheaf cohomology group $H^1(X^d_e,\cF)$ of a 
strip $X^d_e$, where $e<d<b$, is infinite dimensional for any coherent  
analytic sheaf $\cF$. This makes possible the counterexamples of  
Grauert-Andreotti-Rossi. 
On the other hand, the finiteness of this group in dimension greater  
than three implies the Rossi and Andreotti-Siu theorems \cite{Ros:65,AS:70}. 

We shall suppose henceforth without loss of generality that $b>0$.
The function
$\chi=\,-\log{(-\varphi)}$ is
smooth and strictly plurisubharmonic on $X_0$. We set  
\begin{equation} \label{complmetric}
\omega=\,\imat\partial\db\chi=\,-\imat\partial\db\log{(-\varphi)}. 
\end{equation}  
Note that 
$\partial\db\chi=\,\partial\db\varphi/(-\varphi)+ 
(\partial\varphi \wedge\db\varphi)/\varphi^{2}$ and 
$(\partial\varphi \wedge\db\varphi)/\varphi^{2}
=\,\partial\chi\wedge\db\chi$\,. 
Since $\imat{\partial\db\varphi}/{(-\varphi)}$ represents a metric on  
$X_0$, we get the Donnelly-Fefferman condition:
\begin{equation}\label{df} 
\abs{\partial\chi}_{\omega}\leqslant{1}\,. 
\end{equation} 
Since $\chi:X_0\longrightarrow\R$ is proper, \eqref{df} also ensures that  
$\omega$ is complete.  
Indeed, \eqref{df} entails that $\chi$ is Lipschitz  
with respect to the geodesic distance induced by $\omega$, so any geodesic  
ball must be relatively compact.  

We fix in the sequel a regular value  $d\in(-1,0)$ of $\varphi$. The metric $\om$ is complete at  
the pseudoconcave end of $X_d$ and extends smoothly over the boundary 
$bX_d$.  
 
We wish to derive the Poincar\'e inequality for $(0,1)$-forms 
on $X_d$. For this goal we look first at the minus infinity end  
and use the Berndtsson-Siu trick \cite{Be:92,Siu:96}.  
Roughly speaking, it uses the negativity of the trivial line bundle,  
thus avoiding the problems 
raised by the control of the Ricci curvature of $\omega$ at $-\infty$. 
Let us denote by $\cC^{\,0,\,q}_0(X_d)$   
the space of smooth $(0,q)$-forms with compact support in  
$X_d$.  
Let $\vartheta= -*\partial\,{*}$ be the formal adjoint  
of $\db$ with respect to the scalar product  
$(u,v)=\int_{X_d}\langle{u,v}\rangle\,dV_{\omega}$,  
where $\langle{u,v}\rangle=\langle{u,v}\rangle_\om$ and $dV_{\omega}=\omega^n/n!$\,.   
\begin{lemma} \label{gh}
For any $v\in\cC^{\,0,1}_0(X_d)$ we have  
${\norm{v}}^2\leqslant{8}\big({\norm{\db{v}}}^2+{\norm{\vartheta{v}}}^2\big)$.
\end{lemma} 
\begin{proof}   
On the trivial bundle $E=X_d\times\C$ we introduce the auxiliary  
hermitian metric $e^{\chi/2}$. The curvature of $E$ is then
$\Theta(E)=\partial\db{(-\chi/2)}$.  
Let $\vartheta_{\chi}$ be the formal adjoint of  
$\db$ with respect to the scalar product  
$(u,v)_\chi=\,\int_{X_d}\langle{u,v}\rangle\,e^{\chi/2}\,dV_{\omega}$\,. 
Then $\vartheta_{\chi}=e^{-\chi/2}\,\vartheta\,e^{\chi/2}$.  
We apply the Bochner-Kodaira-Nakano formula for $u\in{\cC}^{\,0,1}_0(X_d)$:  
\begin{equation}\label{bkn1} 
\int_{X_d}\big\langle[\imat\,\partial\db{(-\chi/2)},  
\Lambda_{\om}]\,u,u\big\rangle 
\,e^{\chi/2}\,dV_{\omega}\leqslant\,\int_{X_d}\big(|\db{u}|^2+ 
|\vartheta_{\chi}{u}|^2\big)\,e^{\chi/2}\,dV_{\omega}\,, 
\end{equation} 
where $\Lambda_{\om}$ represents the contraction with $\om$ and  
$[A,B]=AB-(-1)^{\deg{A}\cdot\deg{B}}BA$ is the graded commutator  
of the operators $A,B$.  
The idea is to substitute $v=u\,e^{\chi/4}$. 
It is readily seen that 
\begin{equation}\label{dbth1} 
|\db{u}|^2\,e^{\chi/2}\leqslant\,2|\db{v}|^2+\tfrac{1}{8}|\db{\chi}|^2|v|^2\,, \quad
|\vartheta_{\chi}{u}|^2\,e^{\chi/2}\leqslant\,2|\vartheta{v}|^2+\tfrac{1}{8} 
|\partial{\chi}|^2|v|^2 \,. 
\end{equation} 
Moreover 
$
\big\langle[\imat\,\partial\db(-\chi/2),\,\Lambda_{\omega}]u,u\big\rangle 
\,e^{\chi/2}= 
\big\langle[\imat\,\partial\db(-\chi/2),\,\Lambda_{\omega}]v,v\big\rangle\,. 
$
In general, for a $(p,q)$--form $\alpha$ we have the identity  
$\langle[\omega,\Lambda_{\omega}]\alpha,\alpha\rangle=(p+q-n)|\alpha|^2$,
where $n=\dim X$.  
Taking into account that $\omega=\imat\,\partial\db{\chi}$ and that  
$v$ is a $(0,1)$--form, we obtain 
\begin{equation}\label{uv} 
\big\langle[\imat\,\partial\db(-\chi/2),\Lambda_{\omega}]u,u\big\rangle 
\,e^{\chi/2}=\frac{n-1}{2}\,|v|^2\geqslant\frac{1}{2}\,|v|^2\,. 
\end{equation} 
By \eqref{bkn1}, \eqref{dbth1}, \eqref{df}, 
\begin{equation} 
\frac{1}{2}\int_{X_d}|v|^2\,dV_{\omega}\leqslant\,2\int_{X_d}\big(|\db{v}|^2+ 
|\vartheta{v}|^2\big)\,dV_{\omega}+\frac{1}{4}\int_{X_d}|v|^2\,dV_{\omega}\,. 
\end{equation} 
This immediately implies Lemma \ref{gh} for elements $v\in{\cC^{\,0,1}_0(X_d)}$.  
\end{proof}   
 
Let $\eta:(-\infty,0)\longrightarrow\R$ be a smooth 
function such that $\eta(t)=0$ on $(-\infty,-2]$, $\eta'(t)>0$,  
$\eta''(t)>0$  
on $(-2,0)$. Let us introduce the scalar product 
\begin{equation}\label{scalar} 
(u,v)_{\eta(\varphi)}=\int_{X_d}\langle{u,v}\rangle\,e^{-\eta(\varphi)}\,dV_\om\,, 
\end{equation} 
the corresponding norm $\norm{\,\cdot\,}_{\eta(\varphi)}$ and $L^2$ spaces,  
denoted $L^{0,\,q}_2(X_d,\eta(\varphi))$. 
Let $\cC^{\,0,\,q}_0({\overline X}_d)$ be 
the space of smooth $(0,q)$-forms with compact support in ${\overline X}_d$. 

Consider the maximal closed extension of $\db$ to $L^{0,\,q}_2(X_d,\eta(\varphi))$ 
and let $\dbs_{\eta(\varphi)}$ and $\vartheta_{\eta(\varphi)}$ be its Hilbert-space 
and formal adjoints, respectively.    
Then $\vartheta_{\eta(\varphi)}=\vartheta+ 
i\big(\pa\eta(\varphi)\big)$, where $i(\,\cdot\,)$ represents the interior product. 

We denote by $\sigma(P,df)$ the symbol of a differential operator of order one,
calculated on the cotangent vector $df$.  
Then $\sigma(\vartheta,df)=*\,\partial{f}\wedge*$ 
and it is clear that $\sigma(\vartheta_{\eta(\varphi)},df)=\sigma(\vartheta,df)$ 
does not depend on $\eta$. 
Set
$B^{\,0,\,q}= \{\,\alpha\in\cC^{\,0,q}_0(\,{\overline X}_d):\sigma(\vartheta,d 
\varphi)\alpha=0\,\,\text{on} \,\,bX_d\}$.
Integration by parts \cite[Prop. 1.3.1--2]{FoKo1:72} yields 
$\Dom\dbs_{\eta(\varphi)}\cap\cC^{\,0,\,q}_0(\,{\overline X}_d)= 
B^{\,0,\,q}$, $\dbs_{\eta(\varphi)}=\vartheta_{\eta(\varphi)}$
on $B^{\,0,\,q}$. 

\begin{lemma}\label{density} 
The space $B^{\,0,\,q}$ is dense in $\Dom\db\cap\Dom\dbs_{\eta(\varphi)}$ in the graph  
norm 
\begin{equation*}
u\longmapsto\big(\norm{u}_{\eta(\varphi)}^2+\norm{\db{u}}_{\eta(\varphi)}^2+
\norm{\db_{\eta(\varphi)}^*{u}}_{\eta(\varphi)}^2\big)^{1/2}\;.
\end{equation*} 
\end{lemma} 
\begin{proof} 
We use first the idea from \cite[Lemma 4, p.\,92--3]{AV:65} in order to 
reduce the proof to the case of a compactly supported form 
$u\in\cC^{\,0,\,q}_0(\,{\overline X}_d)$.
But then the approximation in the graph norm follows from the  
Friedrichs theorem on the identity of weak and strong derivatives, cf.  
\cite[Prop. 1.2.4]{H1}. 
\end{proof} 
 
We confine next our attention to the fundamental estimate on ${X_d}$.  
\begin{lemma}\label{fdt-est} 
If $\eta$ grows sufficiently fast, there exists a constant $C>0$ such that  
\begin{equation}\label{fund} 
\norm{u}^2_{\eta(\varphi)}\leqslant{C}\Big(\norm{\db{u}}^2_{\eta(\varphi)}+ 
\norm{\dbs_{\eta(\varphi)}{u}}^2_{\eta(\varphi)}+ 
\int_{K}\abs{u}^2\,e^{-\eta(\varphi)}\,dV_\om\Big)\,, 
\end{equation} 
for any $u\in\Dom\db\cap\Dom{\dbs_{\eta(\varphi)}}\subset L^{\,0,\,1}_2(X_d,\eta(\varphi))$,  
where $K=\{-3\leqslant\varphi\leqslant-3/2\}$. 
\end{lemma}  
\begin{proof}
The Morrey-Kohn-H\"ormander estimate 
\cite[Th.\,3.3.5]{H1}, \cite[p.\,429,\,(7.14)]{Gri:66}, 
for $(0,1)$-forms shows that there exists $R>0$ 
such that for sufficiently  
growing $\eta$ the following estimate holds: 
\begin{equation}\label{fund1} 
\norm{u}^2_{\eta(\varphi)}\leqslant R\Big(\norm{\db{u}}^2_{\eta(\varphi)}+\norm{\dbs_{\eta(\varphi)}{u}}^2_{\eta(\varphi)}+ 
\int_{\{-3\leqslant\varphi\leqslant-3/2\}}\abs{u}^2\,e^{-\eta(\varphi)}\,dV_\om\Big)\,, 
\end{equation}  
for $u\in\Dom\db\cap\Dom{\dbs_{\eta(\varphi)}}\subset L^{\,0,\,1}_2(X_d,\eta(\varphi))$, 
$\supp{u}\subset\{-3\leqslant\varphi\}$\,. 
Let $u\in\Dom\db\cap\Dom{\dbs_{\eta(\varphi)}}\subset L^{\,0,\,1}_2(X_d,\eta(\varphi))$. 
The density Lemma \ref{density} shows that to prove \eqref{fund} 
it suffices to consider smooth elements $u$ compactly supported in 
$\overline{X}_d$. 
We choose a cut-off function  
$\rho_1\in\cC^{\infty}(\,\overline{X}_d)$ such that $\supp\rho_1=\{-3\leqslant\varphi\}$, 
$\rho_1=1$ on $\{-2\leqslant\varphi\}$. Set $\rho_2=1-\rho_1$. 
On $\supp\rho_2$, $\eta(\varphi)$ vanishes, therefore $\dbs_{\eta(\varphi)}(\rho_{2}u)= 
\vartheta(\rho_{2}u)$. Upon applying Lemma \ref{gh} for $\rho_{2}u$ we get 
$\norm{\rho_{2}u}^2_{\eta(\varphi)}\leqslant8\big(\norm{\db{(\rho_{2}u)}}^2_{\eta(\varphi)}+ 
\norm{\dbs_{\eta(\varphi)}(\rho_{2}u)}^2_{\eta(\varphi)}\big)$. 
The latter estimate and estimate \eqref{fund1} for $\rho_{1}u$  
together with standard inequalities deliver \eqref{fund}.  
\end{proof}
In the sequel we fix a function $\eta$ as in Lemma \ref{fdt-est}. 
Then the fundamental estimate \eqref{fund} implies the solution 
of the $L^2$ $\db$-Neumann problem.  
Consider the complex of closed, densely defined operators 
\begin{equation*}  
\begin{CD} 
L^{\,0,\,0}_2({X_d},\eta(\varphi))@>{T=\db}>>L^{\,0,\,1}_2({X_d},\eta(\varphi)) 
@>S=\db>>L^{\,0,\,2}_2({X_d},\eta(\varphi))\,, 
\end{CD} 
\end{equation*} 
and the closed, densely defined operator 
\begin{gather*} 
\Dom\lap=\left\lbrace u\in\Dom S\cap\Dom T^*\,:\,Su\in\Dom S^*\,,\; 
T^*u\in\Dom T\right\rbrace\,,\\ 
\lap u=S^*Su+T\,T^*u\quad\text{for}\;u\in\Dom\lap\,. 
\end{gather*}  
We know from a theorem of Gaffney \cite[Prop.\,1.3.8]{FoKo1:72}, that $\lap$ is self-adjoint. 
We denote in the sequel $\cH^{\,0,\,1}=\,\ke{S}\cap\ke{T^*}$. 
\begin{theorem}\label{Neumann}
The following assertions hold true\,{\rm:}
\begin{enumerate} 
\renewcommand{\labelenumi}{(\roman{enumi})} 
\item The operators $T$ and $\lap$ have closed range, $\cH^{\,0,\,1}$
is finite dimensional, and we have the 
strong Hodge decomposition 
\begin{equation*} 
L^{\,0,\,1}_2({X_d},\eta(\varphi))=\Ran(T\,T^*)\oplus 
\Ran(S^*S)\oplus\cH^{\,0,\,1}\,. 
\end{equation*}  
\item There exists a bounded operator $N$ on $L^{\,0,\,1}_2({X_d},\eta(\varphi))$ 
such that $\lap N=N\lap=\Id-P_h$, $P_{h}N=NP_h=0$, where $P_h$ is the  
orthogonal projection on $\cH^{\,0,\,1}$. 
\item If $f\in\Ran T$, the unique solution $u\perp\ke T$ of the equation 
$Tu=f$ is given by $u=T^*Nf$. 
\item The operator $N$ maps  
$L^{\,0,\,1}_2({X_d},\eta(\varphi))\cap\cC^{\,0,\,1}({X_d})$ 
into itself. 
\end{enumerate} 
\end{theorem} 
 
\begin{proof} 
The fundamental estimate \eqref{fund} implies as in \cite[Prop.\,1.2]{Oh:82}
that for any bounded sequence 
$u_k\in\Dom\db\,\cap\,\Dom{\dbs_{\eta(\varphi)}}\subset L^{\,0,\,1}_2(X_d,\eta(\varphi))$
with $\norm{\db u_k}_{\eta(\varphi)}\to0$, $\norm{\dbs u_k}_{\eta(\varphi)}\to0$
one can select a strongly convergent subsequence.
From this follow assertions (i)-(iii) (see e.g.\,\cite[Prop.\,1.1.3]{H1}).
Since $\lap$ is an extension of 
$\db\vartheta_{\eta(\varphi)}+\vartheta_{\eta(\varphi)}\db$, assertion (iv) follows from the 
interior regularity for  
elliptic operators (see e.g.\,\cite[Th.\,2.2.9]{FoKo1:72}).
\end{proof}
\begin{remark}
By using the estimates in local Sobolev norms near the boundary points, 
we can prove as in Folland-Kohn \cite{FoKo1:72} that
$N$ maps $L^{\,0,\,1}_2({X_d},\eta(\varphi))\cap\cC^{\,0,\,1}(\,{\overline{X}_d})$ 
into itself. We could repeat then the solution of the Levi problem as given  
in \cite[Th.\,4.2.1]{FoKo1:72}, in order to find holomorphic  
peak functions, for each boundary point. However, we propose in Corollary
\ref{levi}, a simpler proof for the existence of peak functions, 
which doesn't involve the regularity up to the boundary of the  
$\db$-Neumann problem. 
\end{remark}
\begin{corollary}\label{levi} 
Let $p\in b{X_d}$ and $f$ be a holomorphic function 
on a neigbourhood of $p$ such that $\{f=0\}\cap{\overline{X}_d}=\{p\}$. Then for 
every $m$ big enough, there is a function $g\in 
\cO({X_d})\cap\cC^\infty(\,{\overline{X}_d}\smallsetminus\{p\})$, a smooth function $\Phi$ 
on a neigbourhood $V$ of $p$  and constants 
$a_1,\cdots,a_{m-1}$ such that 
$g=f^{-m}(1+a_{m-1}f+\cdots+ a_1 f^{m-1})+\Phi$
on $V\cap X_d$. In particular, we have $\lim_{z\rightarrow p}\abs{g(z)}=\infty$.
\end{corollary}
\begin{proof} 
Let $U$ be a small neighbourhood of $p$ where $f$ is defined.
Pick $\psi\in\cC^\infty_0(U)$ such that  $\psi=1$ on a neighbourhood $V'$ of $p$.
Set 
$h_m=\psi/f^m$ on $U$, $h_m=0$ on $X\smallsetminus\supp\psi$,
and $v_m=\db h_m$ on $X$.
Observe that $v_m$ belongs to $\cC^{\,0,1}_0(\,\overline{X}_{d+\delta})$  
for $\delta>0$ small enough and $v_m= 0$ on $V'$. Moreover, we have $\db v_m=0$ on $X_{d+\delta}$. 
Fix such a $\delta$ and apply Theorem \ref{Neumann} for $X_{d+\delta}$.  
By this theorem, the 
codimension of $\Ran{T}$ in $\ke{S}$ is finite. For every $m$ big enough, 
there are constants $a_1,\cdots,a_{m-1}$ such that $v=v_m+a_{m-1}v_{m-1}+ 
\cdots + a_1 v_1$ belongs to $\Ran{T}$. Then there is  
$\Phi\in\cC^{\,0,0}(X_{d+\delta})$ such that $\db\Phi=-v$. 
Set $h=h_m+a_{m-1}h_{m-1}+\cdots+a_1{h_1}$ and $g=h+\Phi$. We have 
$\db g=0$ on $X_{d+\delta}\smallsetminus\{f=0\}$. Then $g\in\cO({X_d}) 
\cap\cC^\infty(\,\overline{X}_d\smallsetminus\{p\})$. The function $\Phi$ in the 
corollary is equal to $\Phi$ on $V'$. Thus it is smooth on  
$V:=V'\cap X_{d+\delta}$. The proof is completed.
\end{proof}

\section{The compactification} 
\noindent
In this section we prove Theorem \ref{main}, using the results of Section 
\ref{existence} and the method of \cite{AS:70}. 

\begin{proposition}\label{as} 
Let $d$ be a regular value of $\varphi$. Then for $\delta>0$ small enough 
we have: 
\begin{enumerate} 
\renewcommand{\labelenumi}{(\alph{enumi})} 
\item the holomorphic functions on $X_d$ separate points on $X_{d-\delta}^d$,  
\item the holomorphic functions on $X_d$ give local coordinates  
on $X_{d-\delta}^d$, and
\item for any $e\in(d-\delta,d)$ there exists 
$e^*\in(e,d)$, such that the holomorphically convex hull of 
$X_e$, with respect to the algebra of holomorphic functions on $X_d$, 
is contained in $X_{e^*}$. 
\end{enumerate} 
\end{proposition}

Proposition \ref{as} will be the consequence of the following two lemmas.
We can assume that $bX_{d-\varepsilon}$ is smooth for $\varepsilon>0$  
small enough.
Choose a projection $\pi$ from a neighbourhood of $bX_d$ into $bX_d$.
We will denote by $(x,\varepsilon)$ the point of $bX_{d-\varepsilon}$  
whose projection is $x\in bX_d$.
\begin{lemma}\label{L1}  
Let $x_1$, $x_2$ be two different points in $bX_d$. Then there
are two neighbourhoods $V_1$, $V_2$ of $x_1$, $x_2$ in $bX_d$ and $\nu=\nu(x_1,x_2)>0$ 
such that the holomorphic 
functions of $X_d$ separate $V_1\times (0,\nu]$ and $V_2 \times (0,\nu]$.
\end{lemma}
\begin{proof} 
This is a direct corollary from the existence of a function holomorphic in $X_d$\,,
and ${\cC}^\infty$ in $\overline{X}_d\smallsetminus\{x_1\}$ 
which tends to $\infty$ at $x_1$. 
\end{proof}

\begin{lemma}\label{L2}
Let $x$ be a point of $bX_d$. Then there are a neighbourhood 
$V$ of $x$ in $bX_d$ and $\tau=\tau(x)>0$ 
such that the holomorphic functions in $X_d$ give local coordinates for
$V\times(0,\tau]$. 
\end{lemma}
\begin{proof} 
Without loss of generality and in order to simplify the notations,
we consider the case $n=2$.
Choose a local coordinates system such that $x=0$ and locally $X_d\subset 
\{|z_1-1/2|^2+|z_2|^2< 1/4\}$. 
We now apply Corollary \ref{levi} for functions
$f_1(z)=z_1$ and $f_2(z)=z_1(1-z_2)$.
Denote by $g_1$, $g_2$
the holomorphic functions constructed by this corollary for a number $m$ big enough.
We can also construct the analogue functions if we replace $m$ by $m+1$. Denote by 
 $g'_1$ and $g'_2$ these new functions.
\par
Let $G:X_d\longrightarrow \C^4$ given by $G=(g_1,g_2,g'_1,g'_2)$. 
We will prove that $G$ gives local coordinates. 
Set $I(z)=(z_1 z^{-1}_3,1-z_2 z_3 z^{-1}_1 z^{-1}_4)$.
Let $W$ be a small neigbourhood of $0$.
By Corollary \ref{levi}, the map $I\circ G$ is defined on $W\cap X_d$ and
can be extended to a smooth function on $W$. Moreover, on $W$ 
we have 
$I\circ G(z)=\big(z_1+ O(z_1^2), z_2+O(z_1)\big)$. 
Then $I\circ G$ gives an immersion of $W\cap{X_d}$ in $\C^2$, whenever $W$ is small enough. 
In consequence, $G$ gives coordinates on $W\cap X_d$.
\end{proof}
\begin{proof}[Proof of Proposition \ref{as}]
We cover $bX_d\times bX_d$ by a finite family of 
open sets of the form $V_1\times V_2$ (from Lemma \ref{L1}) and the form 
$V\times V$ (from Lemma \ref{L2}). 
We have a finite family of $\nu$'s and $\tau$'s. Then properties (a) and (b) hold for every 
$\delta$ smaller than these $\tau$'s and $\nu$'s. 
Property (c) is an immediate consequence of Corollary \ref{levi}.
\end{proof}  

\begin{proof}[Proof of Theorem \ref{main}]
First let us remark that the first assertion is a consequence of the second, so we shall prove only the latter.
We assume therefore that the function $\varphi:X\longrightarrow(-\infty,a)$ is strictly plurisubharmonic
everywhere and $a,b>0$.
The proof of the compactification statement for $\dim X\geqslant 3$ in \cite[Prop.\,3.2]{AS:70}
uses only the assertions (a), (b) and (c) of Proposition \ref{as}, so we just have to follow it. 
Namely, let $d$, $\delta$, $e$ and $e^*$ as in Proposition \ref{as}
and denote $P_\varepsilon=\{z\in\C^N\,:\,\abs{z_i}<\varepsilon\}$. 
Proposition \ref{as}
implies as in \cite[Prop.\,3.2]{AS:70} the existence of a holomorphic map 
$\alpha:X_d\to\C^N$ which is an embedding of $X^{e^*}_e$ and 
$\alpha(X_e)\subset P_{1/2}$, $\alpha(\{\varphi=e^*\})\cap P_1=\varnothing$.  
Set $H=\alpha^{-1}(P_1\smallsetminus\overline{P}_{1/2})\cap X_e^{e^*}$.
Since $\alpha(H)$ is a complex submanifold of $P_1\smallsetminus\overline{P}_{1/2}$ of  
dimension at least two, it follows from the Hartogs 
phenomenon \cite[Th. VII,\,D.6]{GuRo:65} that we can find an $\varepsilon\in[1/2,1]$,
such that $\alpha(H)\cap(P_1\smallsetminus\overline{P}_\varepsilon)$ can be extended to 
an analytic subset $V$ of $P_1$.
We can glue the topological spaces $X_d\smallsetminus\alpha^{-1}(\,\overline{P}_\varepsilon)$
and $V$ along $H\smallsetminus\alpha^{-1}(\,\overline{P}_\varepsilon)$ using the  
identification given by the holomorphic map $\alpha$. Hence, we obtain a  
complex space $\widehat{X}_d$\,, such that
$X_d\smallsetminus\alpha^{-1}(\,\overline{P}_\varepsilon)$ and $V$ are open subsets 
of $\widehat{X}_d$\,. 
This turns out to be a Stein space since we can construct a strictly plurisubharmonic exhaustion function, using the function $\varphi$ and the coordinate functions in $\C^N$.
The uniqueness of the Stein completion \cite[Cor. 3.2]{AS:70} entails that $\widehat{X}_d$ does not depend on $d$, so letting $d\to-\infty$ we obtain the desired completion $\widehat{X}$ of $X$.
\end{proof}
\begin{remark} Our method was to embed small strips $X_e^{e^*}$ in $\C^N$  
using holomorphic functions and apply the Hartogs phenomenon.  
One can produce easily holomorphic $(n,0)$-forms on $X_0$ and an embedding 
$\Psi:X_e^{e^*}\to\C\Proj^N$, using the standard $L^2$ estimates 
for $\db$ (cf.\,\cite{Dem:96}). However, it seems that the global Hartogs 
phenomenon in 
$\C\Proj^N$ is an open question \cite[Prob.\,1]{DH97}. Note that by pulling back 
$\Psi(X_e^{e^*})$ to $\C^{N+1}\smallsetminus\{0\}$ we obtain a 
noncompact manifold, so we cannot apply the known results from the 
euclidian space (see also \cite{For:98} where some difficulties of the passing from the local to global Hartogs theorem are exhibited). 
But we can partly use the projective embedding to show Theorem \ref{main}.
Namely, the existence of a non-constant holomorphic function on
$\Psi(X_e^{e^*})\subset\C\Proj^N$ and the arguments from
Sarkis \cite[Cor.\,4.13]{Sa:99} show that we can fill in
$X_e^{e^*}$.
\end{remark}
\begin{remark}[Generalization of Theorem \ref{main}] 
Theorem \ref{main} holds also for normal complex spaces
with isolated singularities (which are the only allowed normal singularities in
dimension two).
Indeed, let $X$ be a hyperconcave end with isolated normal singularities
(Definition \ref{conc} makes sense also for complex spaces). 
Let $\{a_i\}$ denote the singular points and choose functions 
$\varphi_i$ with pairwise disjoint compact supports, such that 
$\varphi_i$ is strictly plurisubharmonic in a neighbourhood of 
$a_i$ and $\lim_{z\rightarrow{a_i}}\varphi_i(z)=-\infty$. Using the function 
$\widehat{\varphi}=\varphi+\sum\varepsilon_i\varphi_i$, with 
$\varepsilon_i$ small enough, we see that $\reg{X}$ 
is a hyperconcave end. 
By Theorem \ref{main} we get a normal Stein completion $Y$
of $\reg{X}$. 
Take $\{V_i\}$ pairwise disjoint Stein neighbourhoods of 
$\{a_i\}$. Then $V_i\smallsetminus\{a_i\}\subset\reg{X}$ and a normal
Stein completion of $V_i\smallsetminus\{a_i\}$ is $V_i$.
Using the uniqueness of a normal Stein completion
\cite[Cor.\,3.2]{AS:70} we infer that the $V_i$ are disjointly
embedded in $Y$. Therefore, $Y$ is also a completion of $X$.
In particular, the singular set of a hyperconcave end with only   
isolated singularities must be finite. 
\end{remark}  
\begin{corollary}\label{antiAS} 
Let $V$ be a Stein manifold, $\dim{V}\geqslant{2}$. Let $K$ be a compact  
completely pluripolar set, $K=\varphi^{-1}(-\infty)$ where $\varphi$ is  
a strictly plurisubharmonic function defined on a neighbourhood $U$ of  
$K$, smooth on $U\smallsetminus{K}$. Then any finite non-ramified covering of $V 
\smallsetminus{K}$ can be compactified to a strongly pseudoconvex space.  
\end{corollary} 
This follows immediately from Theorem \ref{main}, since
$V\smallsetminus{K}$ is a hyperconcave end and any finite non-ramified covering of a  
hyperconcave end is also a hyperconcave end. 
Corollary \ref{antiAS} is in stark contrast to the examples of  
non-compactifiable pseudoconcave ends from \cite{Gra:94,AS:70,Ros:65,Fa:92,Dl93}. 
They are obtained as finite non-ramified coverings of small  
neighbourhoods of the boundaries of some Stein manifolds of dimension two. 
They have `big' holes which cannot be filled, whereas `small', 
i.e. completely pluripolar holes can always be filled.   

\begin{remark}[Complex cobordism]
We can recast Theorem \ref{main} in the light of the cobordism result of 
Epstein and Henkin \cite{EH:01}. Namely,
if $Y$ is a compact strongly pseudoconvex CR manifold of real dimension three,
strictly complex cobordant to $-\infty$, then $Y$ bounds a strongly 
pseudoconvex compact manifold. In particular, $Y$ is embeddable in $\C^N$, for some $N$.  
Note that, if $\dim_{\R} Y>3$ this is automatic by a theorem of Boutet de Monvel 
\cite[p.\,5]{BdM1:74}. On the other side, the examples of Grauert, Andreotti-Siu, Rossi and also Burns \cite{Bu} exhibit compact strongly pseudoconvex CR manifolds of dimension three which do not bound a complex manifold and are not embeddable in $\C^N$.  
\end{remark}  

\section{Compactification by adding finitely many points}

The present section is devoted to proving sufficient conditions for the set
$\widehat{X}\smallsetminus X$ to be analytic. 
In order to prove Theorem \ref{nadel} we consider first the  
particular case when the completion $\widehat{X}$ 
is a Stein space.
 
We begin with some preparations. 
Let $V$ be a complex manifold. We say that $V$ satisfies the  
{\em{Kontinuit\"atssatz}} if 
for any smooth family of closed holomorphic discs $\overline\Delta_t$ in $V$ indexed by 
$t\in [0,1)$ such that $\cup\, b\Delta_t$ lies on a compact subset of $V$, 
 then $\cup\, \overline\Delta_t$ lies on a compact subset of $V$.  
It is clear that every Stein manifold satisfies the Kontinuit\"atssatz, 
using the strictly plurisubharmonic exhaustion function and 
the maximum principle.  
Moreover, if the universal cover of $V$ is Stein
then $V$ satisfies Kontinuit\"atssatz since we can  
lift the family of discs to the universal cover.
\par
Let $F$ be a closed subset of $V$. We say that $F$ is {\em{pseudoconcave}} 
if $V\smallsetminus F$ 
satisfies the local Kontinuit\"atssatz in $V$, i.e. for every $x\in F$ there is a 
neigbourhood $W$ of $x$ such that $W\smallsetminus F$ satisfies the Kontinuit\"atssatz.
Observe that the finite union of pseudoconcave subsets is pseudoconcave and every 
complex hypersurface is pseudoconcave.
 
We have the following proposition which implies the Theorem \ref{nadel}.
\begin{proposition}\label{finite} Let $\widehat\Omega$ be a Stein space  
with isolated singularities $S$ and $K$
a completely pluripolar compact subset of $\widehat\Omega$ which contains $S$.
Assume that 
$\Omega=\widehat{\Omega}\smallsetminus K$ can be 
covered by Zariski-open sets which satisfy the local Kontinuit\"atssatz
in $\widehat{\Omega}\smallsetminus S$.
Then $K$ is a finite set. Moreover, if $\Omega$ itself satisfies the local Kontinuit\"atssatz
in $\widehat{\Omega}\smallsetminus S$, we have $K=S$.
\end{proposition}
\begin{proof} We can suppose
that $\widehat{\Omega}$ is a subvariety of a complex space $\C^N$.
Let $B$ be a ball containing $K$ such that $bB\cap\widehat{\Omega}$ is transversal. 
By hypothesis, we can choose a finite family of Zariski-open sets
$V_1$, $\ldots$, $V_k$ which satisfies the local 
Kontinuit\"atssatz in $\widehat\Omega\smallsetminus S$ and
$\cap\, F_i$ is empty near $bB$, where $F_i={\Omega}\smallsetminus V_i$\,. 
Observe that $F_i$ is 
an analytic subset of $\iom$, $\overline F_i\subset F_i\cup K$. 
Since $F_i\cup (K\smallsetminus S)$ is pseudoconcave in $\widehat\iom\smallsetminus S$,
$F_i$ has no component of codimension $\geqslant2$.
By Hartogs theorem, if $n=\dim{X}>2$, there is a complex 
subvariety $\widehat F_i$ of $\widehat \iom$ which contains $F_i$. 
We will prove this property for the case $n=2$. Set $F=\cup F_i$.
\par
Observe that $\Gamma=F\cap bB$ is an analytic 
real curve. The classical Wermer theorem \cite{We:58} 
says that $\operatorname{hull}(\Gamma)\smallsetminus\Gamma$ is a (possibly void) 
analytic subset of pure dimension $1$ of $\C^N\smallsetminus\Gamma$, where 
$\operatorname{hull}(\Gamma)$ is the polynomial hull of $\Gamma$. 
By the uniqueness theorem, $\operatorname{hull}(\Gamma)\subset \widehat
\iom$. Since $S$ is finite, we have $\operatorname{hull}(\Gamma\cup S)=
\operatorname{hull}(\Gamma)\cup S$.
Set $F'=(F\cup K)\cap\overline B$ and $F''= \operatorname{hull}(\Gamma)\cup S$.
\begin{lemma}[In the case $n=2$]\label{n=2}
 We have $F'\subset F''$.
\end{lemma}
\begin{proof} 
Assume that $F'\not\subset F''$. Then there are a point $p\in F'$ and 
a polynomial $h$ such that
$\sup_{F''}|h|<\sup_{F'}|h|=|h(p)|$. Set $r=h(p)$.
By the maximum principle, we have $h^{-1}(r)\cap F'\subset K\smallsetminus S$. In particular, $p\in K\smallsetminus S$.
Recall that $F'\smallsetminus S$ is pseudoconcave in $\widehat\iom'=
\widehat \iom\cap B\smallsetminus S$.
We will construct a smooth family of discs which does not satisfy
the Kontinuit\"atssatz. This gives a contradiction.
The construction is trivial if $p$ is isolated in $F'$. 
We assume that $p$ is not isolated. By using a small perturbation of $h$,
we can suppose that $h(p)$ is not isolated in $h(F')$.
\par
Set $\Sigma'=h(F')$ and $\Sigma''=h(F'')$.
Then $\Sigma'$ (resp. $\Sigma''$) is included in 
the closed disk (resp. open disk) of center $0$ and of radius $|r|$. 
The holomorphic curves $\{h=\mathrm{const}\}$ define a holomorphic foliation, possibly singular, of $\widehat \iom'$. 
The difficulty is that the fiber $\{h=r\}$ can be singular at $p$.  
Denote by $T$ the set of points $s$ such that $h$ is not 
a submersion on a neighbourhood of $h^{-1}(s)$ on $\widehat\Omega$.
Then $T$ is finite.
\par
Denote also $\Xi$ the unbounded component of $\C\smallsetminus (\Sigma''\cup T)$. It is clear that $\Sigma'$ meets $\Xi$. This property is stable for every small perturbation of the polynomial $h$. Since $K$ is a pluripolar
compact set, $K\cap h^{-1}(a)$ is a polar subset of $h^{-1}(a)$ for every 
$a\in \C$.
\par
Choose a point $b\in\Xi$ such that 
$0<\operatorname{dist}(b,\Sigma')<\operatorname{dist}(b,\Sigma''\cup T)$ 
and $a\in\Sigma'$ such that $\operatorname{dist}(a,b)=\operatorname{dist}(b,\Sigma')$. 
We have $a\not\in \Sigma''\cup T$. Replacing $b$ by a 
point of the interval $(a,b)$ we can suppose that 
$\operatorname{dist}(a,b)<\operatorname{dist}(a',b)$ for every 
$a'\in \Sigma'\smallsetminus\{a\}$. Fix a point $q\in F'$ such that $h(q)=a$. Set $\delta_1=|a-b|$. 
Since $a\not\in T$, we can choose a local coordinates system $(z_1,z_2)$ of an open neigbourhood $W$ of $q$ in $\widehat\Omega^\prime$ such that $z_1=h(z)-b$, $q=(a-b,0)$ and $\{(z_1,z_2),\ |z_1|<\delta_1+\delta_2, |z_2|<2\}\subset W$ with $\delta_2>0$ small enough. 
We can choose a $W$ which does not meet $F''$ and is as small as we want.
\par
Let $L$ be the complex line $\{z_1=a-b\}$.
By the maximum principle, $K'=F'\cap L$ is equal to $K\cap L$.
Then $K'$ is a polar subset of $L$. This implies that
the length of $K'$ is equal to $0$. Thus, for almost every $s\in (0,2)$ the circle $\{|z_2|=s\}\cap L$ does not meet $K'$. Without lost of generality, we can suppose that $K'$ does not meet $\{|z_2|=1\}\cap L$. 
Now we define the disk $\overline\Delta_t$ by
$\overline\Delta_t=\{z_1=(a-b)t, |z_2|\leq1\}$
for $t\in[0,1)$. This smooth familly of discs does not verify the 
local Kontinuit\"atssatz in $W\smallsetminus F'$.
\end{proof}
Now, denote by $\widehat F_i$ the smallest hypersuface of $\widehat\iom$ which contains $F_i$. Set $\widehat F=\cup\,\widehat F_i$.
If $n=2$ we have $F\cup K\subset \widehat F$. This is also true for $n>2$. 
It is sufficient to apply the last lemma for linear slices of $\widehat\iom$.
\begin{lemma}\label{pseudoconcave} 
Let $L$ be a pseudoconcave subset of a complex manifold $V$. If $L$
is included in a hypersurface $L'$ of $V$ then $L$ is itself a hypersurface of $V$.
\end{lemma}
\begin{proof} Observe that $L$ is not included in a subvariety of codimension
$\geqslant2$ of $V$.
Assume that $L$ is not a hypersurface of $V$. Then there is a point $p$ in 
$\operatorname{Reg} L'$ which belongs to the boundary of $L$ in $L'$. 
Choose a local coordinates system 
$(z_1,\ldots,z_n)$ of a neigbourhood $W$ of $p$ such that $W$ contains the unit polydisk
$\Delta^n$, $p\in\Delta^n$ and $L'\cap W=\{z_1=0\}\cap W$. 
We can suppose that $0\not\in L$ and we can choose $W$ small as we want.
\par 
Let $\pi: \Delta^n\longrightarrow \Delta^{n-1}$ be 
the projection on the last $n-1$ coordinates.
Let $q\in L^*=\pi(L\cap \Delta^n)$ such that $\operatorname{dist}(0,L^*)=
\operatorname{dist}(0,q)$. 
Consider the smooth family of discs given by $\overline \Delta_t=\{z=(z_1,z''):|z_1|<1/2\,,z''=tq\}$.
This family does not verify the local Kontinuit\"atssatz in $W\smallsetminus L$.
\end{proof}
\noindent
{\it End of the proof of Proposition \ref{finite}\/.}
We know that $(F_i\cup K)\smallsetminus S$ is pseudoconcave in 
$\widehat\iom\smallsetminus S$ and
$F_i\cup K\subset \widehat F$. By Lemma \ref{pseudoconcave}, 
$(F_i\cup K)\smallsetminus S$ is a hypersurface of $\widehat\iom\smallsetminus S$.
By Remmert-Stein theorem, any analytic set can be extended through a point,
so $F_i\cup K$ is a hypersurface of 
$\widehat\iom$. 
We deduce that $K$ is included in the analytic set
$\cap\,\widehat F_i$ of $\C^N$ which does not intersect $bB$. 
Therefore $K$ must be a finite set. 

If $\Omega$ itself satisfies the local Kontinuit\"atssatz
in $\widehat{\Omega}\smallsetminus S$, we have only one Zariski-open set and the Kontinuit\"atssatz
shows directly that $K=S$.
\end{proof}
\begin{remark} 
The Proposition \ref{finite} holds for $K$ not pluripolar. 
For this case, the proof is more complicated. 
Using another submersion of $\widehat\iom$ in Lemma \ref{n=2}, given by the map 
$z\longmapsto(h(z),h(z)+\varepsilon{z_1},\cdots,h(z)+\varepsilon{z_N})$, 
we can suppose that $R=\max_{K\smallsetminus{S}}\abs{z}> 
\max_{F''}\abs{z}$. Let $q\in{bB_R}\cap(K\smallsetminus{S})$ 
where $B_R$ is the ball of center $0$ and radius $R$. 
Using a small affine change of coordinates, we can suppose that 
$bB_R\cap\widehat\iom$ is transversal at $q$. We then construct easily  
a family of discs close to $T_q(bB_R)\cap\widehat\iom$, which does not 
satisfy the Kontinuit\"atssatz, where $T_q(bB_R)$ is the complex tangent space of  
$bB_R$ at $q$. 
\end{remark} 
\begin{proof}[Proof of Theorem 1.3] 
Let $X$ be a hyperconcave end such that the exhaustion function $\varphi$ 
is overall strictly plurisubharmonic. 
Let $\widehat{X}$ be a smooth completion of $X$. 
Then $\widehat{X}\smallsetminus{X}$ has a strictly pseudoconvex neighbourhood 
$V$. Based on Remmert's reduction theory,  
Grauert \cite[Satz 3,\,p.\,338]{Gra:62} showed that there exists a maximal 
compact analytic set $A$ of $V$. (Note that, by definition, a maximal analytic set has dimension
greater than one at each point.)
Moreover, by \cite[Satz 5,\,p.\,340]{Gra:62} 
there exist a normal Stein space $V'$ with at worst isolated singularities, 
a discrete set $D\subset{V'}$ and a proper holomorphic map  
$\pi:V\longrightarrow{V'}$, biholomorphic 
between $V\smallsetminus{A}$ and with $V'\smallsetminus{D}$ and $\pi(A)=D$. 
That is, $A$ can be blown down to the finite set $D$. 
Of course, $\sing(V')\subset D$.
 
The maximum principle for $\varphi$ implies  
$A\subset\widehat{X}\smallsetminus{X}$. Let $\psi:V'\longrightarrow 
[-\infty,\infty)$ be given by $\psi=\varphi\circ\pi^{-1}$ on $V'\smallsetminus{D}$  
and $\psi=-\infty$ on $\pi(\widehat{X}\smallsetminus{X})$. 
Then $\psi$ is a strictly plurisubharmonic function on $V'$ and 
$\pi(\widehat{X}\smallsetminus{X})$ is its pluripolar set. 
By Proposition \ref{finite}, 
$\pi(\widehat{X}\smallsetminus{X})$ is a finite set. 
Therefore $\widehat{X}\smallsetminus{X}$ consists of $A$ and 
possibly a finite set $D'$. If $X$ has a Stein cover, it follows from the
Kontinuit\"atssatz that $\pi(\widehat{X}\smallsetminus{X})=\sing(V')$. 
Therefore $D'=\varnothing$ and $D=\sing(V')$.
\end{proof}

\begin{remark}
If in Theorem \ref{nadel} we suppose only that $X$ admits a Zariski-open 
dense set which is uniformized by a Stein manifold, we can prove in the same way,
that $\widehat{X}\smallsetminus{X}$ is included in a hypersurface of 
$\widehat{X}$, i.e. $X$ contains a Zariski-open dense set of $\widehat{X}$.
\end{remark}
\section{Extension of Nadel's theorems}
Our goal is to extend to dimension two Nadel's theorems \cite{Na:90}.
The arithmetic quotients are, with a few exceptions, pseudoconcave manifolds
carrying a positive line bundle. In this respect, Nadel \cite{Na:90}
considered the realization as a quasiprojective manifold of a class
of manifolds $X$ with hyperconcave ends and of dimension greater than three.
The method of Nadel is to compactify the manifold
by the theorem of Rossi, and then to apply differential geometric methods,
like the existence of K\"ahler-Einstein metric and the Schwarz-Pick 
lemma of Yau and Mok-Yau. 
Corollary \ref{nad2} is a generalization of
\cite[Th. 0.2]{Na:90}.
\begin{proof}[Proof of Corollary \ref{nad2}]
The necessity of conditions (i) and (ii) is obvious, while the necessity of (iii)
follows from a theorem of Griffiths \cite[Th. I]{Gri:71}.

For the sufficiency, we need the embedding theorem of 
Andreotti-Tomassini \cite[Th. 2, p.\,97]{AnTo:70}, \cite[Lemma 2.1]{NT:88}  
(see also \cite[Th. 4.1]{AS:70}). 
This theorem shows that condition (ii) implies the embedding of $X$ as an open set in a 
smooth projective manifold $\widehat{X}$.
We conclude by Theorem \ref{nadel}.
\end{proof}

Our result pertains to the work of Nadel and Tsuji \cite{NT:88} which generalizes the 
compactification of arithmetic quotients of any rank, by showing that certain 
pseudoconcave manifolds are quasi-projective. In dimension two 
their condition coincides with hyperconcavity. 
Corollary \ref{nad2} yields an extension of their theorem
in dimension two, together with a completely complex-analytic proof
of the compactification of arithmetic quotients. 

As a consequence of Corollary \ref{nad2} we also get a slightly stronger form of 
\cite[Main Theorem]{SY:82} (also noted by Nadel in dimension grater than
three):
\begin{corollary}\label{SiuYau}
Let $X$ be a complete K\"ahler manifold of finite volume and bounded 
negative sectional curvature.
If $\dim X\geqslant 2$, $X$ is biholomorphic to a quasiprojective 
manifold which can be compactified to a Moishezon space by adding finitely many singular points.
\end{corollary}
\begin{proof}
Indeed, the same argument as in \cite{SY:82} or  
\cite[\S3]{NT:88} shows, with the help of the Busemann function, 
that $X$ is hyper $1$-concave. Moreover, the negativity of the curvature  
implies that the canonical bundle $K_X$ is positive and the universal cover of $X$ is Stein. 

We show that the positivity of $K_X$ implies the ring $\oplus_{k>0}H^0(X,K^k_X)$ separates points and 
gives local coordinates everywhere on $X$. 
If $\varphi$ denotes the exhaustion function of $X$, 
$\imat\big(A\Theta(K_X)+\partial\db(-\log(-\varphi))\big)$, $A\gg1$, is a complete  
K\"ahler metric on $X$, where $\Theta(K_X)$ is the curvature of $K_X$.
By the $L^2$ estimates with singular weights for positive line bundles 
(see e.g. Demailly \cite{Dem:96}) our contention follows.
We can thus apply Corollary \ref{nad2}.
\end{proof}

In the same vein as in Corollary \ref{nad2}, we can show that Nadel's main result 
\cite[Th.\,0.1]{Na:90} holds in dimension two. 
Namely, let $X$ be a connected Moishezon, hyper 1-concave manifold of dimension $n\geqslant 2$,
which can be covered by Zariski-open sets uniformized by Stein manifolds.
Then $X$ can be compactified by adding finitely many points to a compact
Moishezon space.

We close the section with one more result  about embedding of  
hyper $1$-concave manifolds.  
Generalizing the Andreotti-Tomassini theorem, Andreotti-Siu  
\cite[Th.\,7.1]{AS:70} show that a strongly $1$-\,concave  
manifold $X$ of $\dim{X}\geqslant{3}$ can be embedded in the projective space,
if it admits a line bundle $E$ such that  
$\oplus_{k>0}{H^0}(X, E^k)$ gives local coordinates on a sufficiently large  
compact of $X$. The proof is  
based on techniques of extending analytic sheaves. Moreover, they show through
an example \cite[p.\,267--70]{AS:70} that the result  
breaks down in dimension two. 
We prove now however, that if we impose the condition of hyperconcavity,  
the result occurs also in dimension two. 
Here $\varphi$ and $b$ have the same meaning as in Definition 
\ref{conc}.
\begin{proposition}\label{hyper-as71}
Let $X$ be a hyper $1$-concave manifold of dimension $n\geqslant 2$. 
Let $c$ be a real number 
such that $c<b$. Assume there is a line bundle $E$ over $X'=\{\varphi>c\}$ such that the 
ring $\oplus_{k>0} H^0(X',E^k)$ gives local coordinates on $X'$. Then $X$ is biholomorphic to an
open subset of a projective manifold.
\end{proposition}

\begin{proof} By Theorem \ref{main}, $X$ is an open subset of a variety $\widehat X$ with isolated 
singularities. Moreover $\widehat X\smallsetminus X'$ is a Stein space.
Replacing $c$ by a $c'$ such that $c<c'<b$ we can suppose that there are holomorphic 
sections $s_0$, $\ldots$, $s_m$ of $H^0(X',E^k)$ which give local coordinates of $X'$ where 
$k$ is big enough. We can define a holomorphic map $\pi:X'\longrightarrow \Proj^m$ 
by 
$\pi(z)=[s_0(z):\cdots:s_m(z)]$.
Then $\pi$ gives a local immersion of $X'$ in $\Proj^m$.
Since $\widehat X\smallsetminus X'$ is embeddable in an euclidian space, a theorem of 
Dolbeault-Henkin-Sarkis \cite{DH97,Sa99} implies that $\pi$ can be extended to a meromorphic 
map from $\widehat X$ into $\Proj^m$.
\par
Denote by $Z$ the set consisting of the singular points of $\widehat X$,  
the points of indeterminacy of $\pi$ and the 
critical points of $\pi$. Then $Z$ is a compact analytic subset of $\widehat X\smallsetminus X'$. 
Since $\widehat X\smallsetminus X'$ is Stein space, $Z$ is a finite set. 
The map $\pi$ gives local immersion of $\widehat X\smallsetminus Z$ in $\Proj^m$.
Let $H$ be the hyperplane line bundle of $\Proj^m$ and set $L=\pi^*(H)$. Then $L$ is a positive line bundle of $\widehat X\smallsetminus Z$. In particular $L$ is positive on $X\smallsetminus Z$ and by a theorem of Shiffman \cite{Sh:71} extends to a positive line bundle on $X$.
By the argument in the proof of Corollary \ref{SiuYau},
we show that $\oplus_{k>0} H^0(X,L^k\otimes K_X)$ separates points and gives local coordinates on $X$. By \cite[Lemma 2.1]{NT:88} $X$ is biholomorphic to an open subset of a projective manifold.
\end{proof}

\subsection*{Acknowledgements}
We thank Professors Mihnea Col\c{t}oiu and
Bo Berndtsson for discussions and constant encouragement
and Professor J\"urgen Leiterer for
hospitality during the preparation of this paper. 
We are grateful the referee for helpful suggestions.


\providecommand{\bysame}{\leavevmode\hbox to3em{\hrulefill}\thinspace}
\providecommand{\MR}{\relax\ifhmode\unskip\space\fi MR }
\providecommand{\MRhref}[2]{%
  \href{http://www.ams.org/mathscinet-getitem?mr=#1}{#2}
}
\providecommand{\href}[2]{#2}

\end{document}